\documentclass[11.05pt,a4paper]{amsart}
\usepackage{amsthm,amsfonts,amssymb,mathrsfs,amsmath}
\usepackage{latexsym}
\usepackage{xspace}
\usepackage[all]{xy}
\usepackage{longtable}
\usepackage{amscd}

\headsep=1cm \numberwithin{equation}{section}
\newtheorem{theorem}{Theorem}[section]
\newtheorem{lemma}[theorem]{Lemma}

\newtheorem{corollary}[theorem]{Corollary}

\theoremstyle{definition}
\newtheorem{definition}[theorem]{Definition}
\theoremstyle{remark}
\newtheorem{remark}[theorem]{Remark}
\newtheorem{example}[theorem]{Example}
\newtheorem{question}[theorem]{Question}

\newcommand{\Ass}{\operatorname{Ass}}
\newcommand{\im}{\operatorname{im}}
\newcommand{\grade}{\operatorname{grade}}

\newcommand{\Assh}{\operatorname{Assh}}
\newcommand{\Spec}{\operatorname{Spec}}

\newcommand{\rank}{\operatorname{rank}}

\newcommand{\pd}{\operatorname{pd}}

\newcommand{\Gdim}{\operatorname{G--dim}}
\newcommand{\CIdim}{\operatorname{CI--dim}}

\newcommand{\Ho}{\operatorname{H}}

\newcommand{\Syz}{\operatorname{Syz}}

\newcommand{\G}{\operatorname{G}}

\newcommand{\Ext}{\operatorname{Ext}}

\newcommand{\Tor}{\operatorname{Tor}}
\newcommand{\Hom}{\operatorname{Hom}}

\newcommand{\depth}{\operatorname{depth}}

\newcommand{\coker}{\operatorname{coker}}

\newcommand{\lo}{\longrightarrow}
\newcommand{\fm}{\frak{m}}
\newcommand{\fp}{\frak{p}}

\newcommand{\fn}{\frak{n}}

\newenvironment{prf}[1][Proof]{\begin{proof}[\bf #1]}{\end{proof}}

\begin{document}

\author[K. Divaani-Aazar, F. Mohammadi Aghjeh Mashhad, E. Tavanfar and M. Tousi]
{Kamran Divaani-Aazar, Fatemeh Mohammadi Aghjeh Mashhad, Ehsan Tavanfar\\ and\\ Massoud Tousi}

\title[On the New Intersection Theorem for totally reflexive modules]
{On the New Intersection Theorem for totally reflexive modules}

\address{K. Divaani-Aazar, Department of Mathematics, Alzahra University, Vanak, Post Code
19834, Tehran, Iran-and-School of Mathematics, Institute for Research in Fundamental Sciences
(IPM), P.O. Box 19395-5746, Tehran, Iran.}
\email{kdivaani@ipm.ir}

\address{F. Mohammadi Aghjeh Mashhad,  Parand Branch, Islamic Azad University, Tehran, Iran.}
\email{mohammadi\_fh@yahoo.com}

\address{E. Tavanfar, School of Mathematics, Institute for Research in Fundamental Sciences (IPM),
P.O. Box 19395-5746, Tehran, Iran.}
\email{tavanfar@ipm.ir}

\address{M. Tousi, Department of Mathematics, Shahid Beheshti University, G.C., Tehran,
Iran-and-School of Mathematics, Institute for Research in Fundamental Sciences (IPM), P.O. Box
19395-5746, Tehran, Iran.}
\email{mtousi@ipm.ir}


\keywords {Big Cohen-Macaulay module; complete intersection dimension; deformation; Gorenstein dimension;
quasi-deformation; totally reflexive module.\\
The first and fourth authors were supported by grants from IPM (no. 92130212 and no. 90130211; respectively).
The research of third author is also supported by IPM}

\begin{abstract} Let $(R,\fm,k)$ be a local ring. We establish a totally reflexive analogue of the New
Intersection Theorem, provided for every totally reflexive $R$-module $M$, there is a big Cohen-Macaulay
$R$-module $B_M$ such that the socle of $B_M\otimes_RM$ is zero. When $R$ is a quasi-specialization of a
$\G$-regular local ring or when $M$ has complete intersection dimension zero, we show the existence of such
a big Cohen-Macaulay $R$-module. It is conjectured that if $R$ admits a non-zero Cohen-Macaulay module of
finite Gorenstein dimension, then it is Cohen-Macaulay. We prove this conjecture if either $R$ is a
quasi-specialization of a $\G$-regular local ring or a quasi-Buchsbaum local ring.
\end{abstract}

\maketitle

\section{Introduction}

Throughout, $(R,\fm,k)$ is a commutative Noetherian local ring with identity. The celebrated New Intersection
Theorem is perceived as a deep
result at the interface of homological and local algebra.  It provides simple proofs for several outstanding
homological conjectures e.g. Auslander's zero-divisor conjecture \cite[Theorem 6.2.3]{Ro} and Bass' conjecture
\cite[Theorem 5.1]{PS}. The New Intersection Theorem was proved in prime characteristic by Peskine and Szpiro
\cite{PS} in 1973. Then Hochster's works \cite{H1} and \cite{H2} established a reduction to prime characteristic
from equicharacteristic zero to give a proof of this theorem in every equicharacteristic ring in 1975. Finally,
in 1987, Roberts \cite{Ro} proved the New Intersection Theorem for mixed characteristic rings by using local
Chern characters.

The New Intersection Theorem asserts that if $$0\lo F_s\lo F_{s-1}\lo \cdots \lo F_1\lo F_0\lo 0$$
is a non-exact complex of finitely generated free $R$-modules with finite length homology modules, then
$\dim R\leq s$. Using the New Intersection Theorem, one can easily see that if $R$ admits a non-zero Cohen-
Macaulay module with $\pd_RM<\infty$, then $R$ must be Cohen-Macaulay.

One of the most important notions in Gorenstein homological algebra is that of totally reflexive modules, which
was defined by Auslander \cite{Au}. The local ring $(R,\fm,k)$ is Gorenstein if and only if $k$ admits a finite
resolution by totally reflexive $R$-modules.  On the other hand, over a Gorenstein local ring, the totally
reflexive modules are precisely the maximal Cohen-Macaulay modules. From this perspective, over a Gorenstein
local ring, a totally reflexive module is regarded as a common generalization of a finitely generated free module
and a maximal Cohen-Macaulay module. Since many results in classical homological algebra have counterparts
in Gorenstein homological algebra, the following questions are raised naturally:

\begin{question}\label{1.1} Let $(R,\fm,k)$ be a local ring. Assume that $$0\lo G_s\lo G_{s-1}\lo \cdots\lo G_1\lo
G_0\lo 0$$ is a non-exact complex of totally reflexive $R$-modules with finite length homology modules.
Is $\dim R\leq s?$
\end{question}

\begin{question}\label{1.2} Let $(R,\fm,k)$ be a local ring. Assume that $R$ admits a non-zero Cohen-Macaulay module
of finite Gorenstein dimension. Is $R$ Cohen-Macaulay?
\end{question}

Question 1.2 was asked in \cite[page 40]{C1}, \cite[Question 1.31]{CFH} and \cite{T2}.
This question has been studied by many authors; see  e.g. \cite{T3}, \cite{T2}, \cite{GHT} and \cite{DMT}.

This paper is concerned with the study of these questions. Note that Question 1.1 easily implies Question 1.2.
Theorems \ref{2.4}, \ref{2.7} and \ref{3.6} are the main results of this paper. Let $$G=0\lo G_s\lo G_{s-1}\lo
\cdots\lo G_1\lo G_0\lo 0$$ be a non-exact complex of totally reflexive $R$-modules with finite length homologies.
In Theorem \ref{2.4}, we show that $\dim R\leq s$, provided for every totally reflexive $R$-module $M$, there is
a big Cohen-Macaulay $R$-module $B_M$ such that the socle of $B_M\otimes_RM$ is zero.  We show the existence of such
a big Cohen-Macaulay $R$-module when $R$ is a quasi-specialization of a $\G$-regular local ring or when $M$
has complete intersection dimension zero; see Theorem \ref{2.7}. Also in Corollary \ref{2.8}, we show that if $\CIdim_RG
<\infty$, then $\dim R\leq s$.  Note that this generalizes the New Intersection Theorem without imposing any extra
assumption on the ring $R$.  The notion of $\G$-regular rings was introduced by Ryo Takahashi in \cite{T1}. In Remark
\ref{2.10}, we provide many examples of quasi-specializations of $\G$-regular rings which are neither $\G$-regular nor
Cohen-Macaulay. We provide an affirmative answer to Question 1.2 in the cases $R$ is a quasi-specialization of a
$\G$-regular local ring or a quasi-Buchsbaum local ring; see Corollary \ref{3.1} and Theorem \ref{3.6}.

\section{Question 1.1}

In this section, we are dealing with Question \ref{1.1}. Let us begin by recalling some needed definitions.

For an $R$-complex $X$,  we set $\inf X:=\inf \{i\in \mathbb{Z}\ |\ \Ho_i(X)\neq 0\}$ and $\sup X:=\sup \{i\in
\mathbb{Z}\ |\ \Ho_i(X)\neq 0\}$. We obey the convention that the infimum and supremum of the empty set are
$+\infty$ and $-\infty$; respectively. In the sequel, we use $``\simeq"$ for denoting a quasi-isomorphism between
two complexes. An $R$-complex $X$ is said to be \emph{homologically finite} if every homology module of $X$ is
finitely generated and $\Ho_i(X)=0$ for every $|i|\gg 0$. Finally, {\it dimension} of an $R$-complex $X$ is defined
as: $$\dim_RX:=\sup \{\dim_R\Ho_i(X)-i|\ i\in \mathbb{Z} \}.$$

\begin{definition}\label{2.1}
\begin{enumerate}
\item[i)] A finitely generated $R$-module $M$ is said to be \emph{totally reflexive} if there exists
an exact complex $$F=\cdots\rightarrow F_2\overset{d_2}{\rightarrow}F_1\overset{d_1}{\rightarrow}
F_0\overset{d_0}{\rightarrow} F_{-1}\rightarrow \cdots,$$ of finitely generated free $R$-modules
such that $\Hom_R(F,R)$ is exact and $M\cong \im\ d_0$. Such an $R$-complex $F$ is called a \emph{complete
free resolution} of $M$. For each integer $i$, we set $\Syz_i^F(M):
=\im d_i$.
\item[ii)] For a homologically finite $R$-complex $M$, the {\it Gorenstein dimension} of $M$, $\Gdim_RM$,
is defined as the
infimum of all integers $n$ such that there exists a complex  $$0\lo G_n\lo G_{n-1}\lo \cdots \lo
G_{l+1}\lo G_{l}\lo 0$$ of totally reflexive $R$-modules with $G\simeq M$ and $G_n\neq 0$.
\item[iii)] The class of all totally reflexive $R$-modules is denoted by $\G(R)$. The ring $R$ is said
to be \text{G}-\emph{regular} if every $R$-module $M\in \G(R)$ is free.
\end{enumerate}
\end{definition}

Every finitely generated projective module is totally reflexive. So for a finitely generated
$R$-module $M$, one has $\Gdim_RM\leq \pd_RM$, moreover, equality holds if either $R$ is $\G$-regular or
$\pd_RM<\infty$; see \cite[Proposition 1.8 (2)]{T1} and \cite[Proposition 1.2.10]{C1}.

Recall that an $R$-module $B$ is called {\em big Cohen-Macaulay} if there exists a system of parameters
$\mathbf{x}=x_{1},\ldots, x_{n}$ for $R$ such that $\mathbf{x}$ is a $B$-regular sequence.  If every system of
parameters for $R$ is a $B$-regular sequence, then $B$ is called {\em balanced big Cohen-Macaulay}. Recently,
Andr\'{e} \cite{An} shows that any local ring $R$ possesses a big Cohen-Macaulay module which is an $R$-algebra.
Note that, by \cite[Corollary 8.5.3]{BH}, if $B$ is a big Cohen-Macaulay $R$-module, then $\widehat{B}$, the
$\fm$-adic completion of $B$, is a balanced big Cohen-Macaulay $R$-module.

\begin{lemma}\label{2.2}  Let $(R,\fm,k)$ be a local ring and $B$ a big Cohen-Macaulay $R$-module.
Then $B\otimes_RM\neq 0$ for every non-zero finitely generated $R$-module $M$.
\end{lemma}

\begin{prf} Set $d:=\dim R$ and let $\mathbf{x}=x_1,\ldots, x_d\in \fm$ be a system of parameters for $R$ such
that $\mathbf{x}$ is a $B$-regular sequence. Thus $B/(\mathbf{x})B \neq 0$ and then we can easily verify that
$B/\fm B\neq 0.$  Hence $$k\otimes_R(B\otimes_R M)\cong B/\fm B\otimes_k M/\fm M\neq 0,$$
for every non-zero finitely generated $R$-module $M$.
\end{prf}

Recall that for a not necessarily finitely generated $R$-module $M$, the {\it depth} of $M$ is defined by
$$\depth_RM:=\inf \{i\in\mathbb{N}_0\ |\ \Ext_R^i(k,M)\neq 0\}.$$ (Note that by our convention on the infimum
of the empty set, we have $\depth_R0=+\infty$.) By \cite[Theorem 6.1.6]{St}, this definition coincides with
\cite[Definition 9.1.1]{BH}.

Next, we recall the Peskine-Szpiro Acyclicity Lemma.

\begin{lemma}\label{2.3} Let $(R,\fm,k)$ be a local ring. Let  $$Y:=0\lo Y_{s}\lo Y_{s-1}\lo \cdots \lo
Y_1\lo Y_0$$ be a complex of $R$-modules such that for each $i=1,\ldots, s$, one has $\depth Y_i\geq i$
and either $\Ho_i(Y)=0$ or $\depth_R(\Ho_i(Y))=0$. Then $\Ho_i(Y)=0$ for each $i\geq 1$.
\end{lemma}

\begin{prf} See e.g. \cite[Proposition 1.1.1]{St}.
\end{prf}

For an $R$-module $M$, its {\it socle} is defined as $\text{Soc}(M):=(0:_M\fm).$

Let $X$ be an $R$-complex such that the $R$-module $\Ho_i(X)$ is finitely generated for every $i$ and $\Ho_i(X)=0$
for every $i\ll 0$. Set $t:=\inf X$. By \cite[Theorem A.3.2 (L)]{C1}, there is a complex $$F:=\cdots \lo F_i\lo F_{i-1}\lo
\cdots \lo F_{t+1}\lo F_t\lo 0$$ of finitely generated free $R$-modules such that $F\simeq X$. A such complex $F$
is called a free resolution of $X$.

We are ready to present our first main result in this section.

\begin{theorem}\label{2.4} Let $(R,\fm,k)$ be a local ring. Assume that to each totally reflexive $R$-module $G$, a big
Cohen-Macaulay $R$-module $B_{G}$ is assigned such that $\text{Soc}(B_G\otimes_RG)=0$. Then existences of a non-exact
complex $$\mathcal{G}:=0\rightarrow G_{s}\rightarrow G_{s-1}\rightarrow\cdots\rightarrow G_{0}\rightarrow0,$$ of finite
length homologies and consisting of totally reflexive $R$-modules, yields $s\ge \dim R$.
\end{theorem}

\begin{prf} By considering the free resolution of the complex $\mathcal{G}$ and applying \cite[Theorem 2.3.7]{C1}, we may reduce to
the special case that $$\mathcal{G}:=0\rightarrow G\rightarrow F_{s-1}\rightarrow F_{s-2}\rightarrow\cdots\rightarrow
F_{0}\rightarrow 0$$ for some finite free modules $F_{0},\ldots,F_{s-1}$, as well as a totally reflexive module $G$.
Let $\text{Spec}(R)^{\circ}$ denote, the punctured spectrum of $R$, $\text{Spec}(R)\setminus \{\fm\}$ and $(-)^*$ denote the
$R$-dual functor $\Hom_R(-,R)$.

Since all homology modules of $\mathcal{G}$ have finite length, it turns out that $\mathcal{G}$ is locally split-exact on
$\text{Spec}(R)^{\circ}$, and so $G$ is locally free on $\text{Spec}(R)^{\circ}$.
As $\mathcal{G}$ is locally split-exact on $\text{Spec}(R)^{\circ}$, it follows that $\mathcal{G}^{*}$ is also locally
split-exact on $\text{Spec}(R)^{\circ}$, and so all of its homology modules have finite length.
If $\mathcal{G}^{*}$ is exact, then we would have $\pd_{R}(G^{*})<\infty$, so $G^{*}$ and thus $G$ would be free implying
that $s\ge \dim R$, by the miracle of the New Intersection Theorem. Hence, it suffices to consider only the case where
$\mathcal{G}^{*}$
is non-exact. If necessary, we may replace $\mathcal{G}^{*}$ with the complex
$$0\rightarrow F_{0}^{*}\rightarrow\cdots\rightarrow F_{s-3}^{*}\rightarrow F_{s-2}^{*}\rightarrow \ker(F_{s-1}\rightarrow
G^{*})\rightarrow0,$$ and so on, to impose the condition $F_{s-1}^{*}\rightarrow G^{*}$	has non-zero finite length cokernel.
Pick a big-Cohen-Macaulay module $R$-module $B:=B_{G^{*}}$ such that $\text{Soc}(B\otimes_RG^{*})=0$. Lemma 2.2 implies that
$$\Ho_{0}(B\otimes_R\mathcal{G}^{*})=B\otimes_R\Ho_{0}(\mathcal{G}^{*})\neq 0,$$ and so the complex $B\otimes_R\mathcal{G}^{*}$
is non-exact. It has the form $$\mathcal{C_{\bullet}=}0\rightarrow C_{s+1}:=B^{n_{0}}\rightarrow
C_{s}:=B^{n_{1}}\rightarrow\cdots\rightarrow C_{2}=B^{n_{s-1}} \rightarrow C_{1}=B\otimes_RG^{*}\rightarrow C_{0}=0.$$
Since the complex $\mathcal{C_{\bullet}}$ is locally split-exact on $\text{Spec}(R)^{\circ}$, all of its homology modules
are either zero or of depth zero. In case, $s\le\text{dim}\ R-1$, then we have $\depth_RC_{i}\ge i$ for every
$i\ge1$, so we may apply Lemma \ref{2.3} to conclude that the complex $C_{\bullet}$ is exact which is a contradiction.
\end{prf}

Now, we provide a sufficient condition for satisfying the assumption of Theorem \ref{2.4}. Also, it will be used in the
proof of Theorem \ref{2.7}.

\begin{lemma}\label{2.5} Let $(R,\fm,k)$ be a local ring, $M$ a totally reflexive $R$-module and $\ell$ an integer with
$1\leq \ell \leq \dim R$. Let $\Syz_i^F(M)$ be as in Definition \ref{2.1}. Assume that there exists a big Cohen-Macaulay
$R$-module $B_M$ such that $\Tor^{R}_{1}(B_M,\Syz_{-i}^F(M))=0$ for all $1\le i\le \ell$. Then $\depth_R(B_M\otimes_RM)
\geq \ell$.
\end{lemma}

\begin{prf} Set $d:=\dim R$ and $B:=B_M$.  By \cite[Exercise 9.1.12 (a),(b)]{BH} for every $R$-module $X$ with $X\neq \fm X$,
one has $\depth_RX\leq d$ and $\depth_R0=+\infty$. Also, by \cite[Proposition 9.1.2 (e)]{BH} every short exact sequence
$$0\lo X\lo Y\lo Z\lo 0,$$ yields an inequality $$\depth_RX\geq  \min \{\depth_RY, \depth_RZ+1\}. \  \ (\dag)$$

Since $M$ is totally reflexive, there exists an exact complex $$F=\cdots\rightarrow F_2\overset{d_2}{\rightarrow}
F_1\overset{d_1}{\rightarrow} F_0\overset{d_0}{\rightarrow}F_{-1}\rightarrow \cdots,$$ of finitely generated free
modules such that $\Hom_R(F,R)$ is exact and $M\cong \im\ d_0$. By the assumption, we have the short exact sequences
$$0\lo B\otimes_R\Syz^F_{1-i}(M)\lo B\otimes_RF_{-i}\lo B\otimes_R\Syz^F_{-i}(M)\lo 0$$ for all $1\leq i\leq \ell$. Note
that $B\otimes_R\Syz^F_{0}(M)\cong B\otimes_RM$. We may and do assume that $M\neq 0$.

First, assume that $$\depth_R(B\otimes_R\Syz^F_{-i}(M))+1<d$$ for all $1\leq i\leq \ell$. Since $\depth_RB=\dim R$,  applying
$(\dag)$ to the short exact sequences $$0\lo B\otimes_R\Syz^F_{1-i}(M)\lo B\otimes_RF_{-i}\lo B\otimes_R\Syz^F_{-i}(M)\lo 0;
\
\   i=1,\dots, \ell,$$ successively yields that
\begin{align*}
\depth_R(B\otimes_RM)  &\geq  \depth_R(B\otimes_R\Syz^F_{-1}(M))+1 \\
&\geq \depth_R(B\otimes_R\Syz^F_{-2}(M))+1+1 \\
&\vdots \\
&\geq \depth_R(B\otimes_R\Syz^F_{-\ell}(M))+1+\ell-1\\
&\geq \ell.
\end{align*}

Next, assume that there is an integer $1\leq j\leq \ell$ such that $$d\leq \depth_R(B\otimes_R\Syz^F_{-j}(M))+1, $$ and set
$$s:=\inf\{i\in \{1,2, \dots, \ell\}\ |\ d\leq \depth_R(B\otimes_R\Syz^F_{-i}(M))+1\}.$$ By repeating the above inequalities
$s$ times,
we get $$\depth_R(B\otimes_RM)\geq d+s-1\geq \ell+s-1\geq \ell.$$
\end{prf}

Next, we recall the definition of complete intersection dimension; see \cite{AGP} and \cite{Sa}.

\begin{definition}\label{2.6} Let $(R,\fm,k)$ be a local ring.
\begin{enumerate}
\item[i)] A \emph{deformation} is a surjective homomorphism of local rings $R\leftarrow A$ with
the kernel generated
by an $A$-sequence. In this situation, we say $R$ is a \emph{specialization} of $A$.
\item[ii)] A \emph{quasi-deformation} is a diagram of local homomorphisms $R\rightarrow R'\leftarrow A$ such that
the first map is flat and the second map is a deformation. In this situation, we say $R$ is a \emph{quasi-specialization} of $A$.
\item[iii)] For a homologically finite $R$-complex $M$, \emph{complete intersection dimension}
of $M$, $\CIdim_RM$, is defined as
$$\CIdim_RM:=\inf\{\pd_A(R'\otimes_RM)-\pd_AR' \ | \ R\rightarrow R'\leftarrow A \ \text{is a
quasi-deformation} \}.$$
\end{enumerate}
\end{definition}

For a homologically finite $R$-complex $M$, \cite[Proposition 3.3]{Sa} implies that $$\Gdim_RM\leq \CIdim_RM\leq \pd_RM$$
and if one of these dimensions is finite, then it equals those to its left.

The next is our second main result in this section.

\begin{theorem}\label{2.7} Let $(R,\fm,k)$ be a local ring.
\begin{enumerate}
\item[i)] Assume that $M$ is a finitely generated $R$-module with $\CIdim_RM=0$. Then there is a balanced big Cohen-Macaulay
$R$-module $B_M$ such that $\depth_R(B_M\otimes_RM)\geq \dim R$.
\item[ii)] Assume that $A$ is a $\G$-regular local ring and there is a quasi-deformation $R\rightarrow R'\leftarrow A$. Then,
there is a balanced big Cohen-Macaulay $R$-module $B$ such that $\depth_R(B\otimes_RM)\geq \dim R$ for every totally reflexive
$R$-module $M$.
\end{enumerate}
\end{theorem}

\begin{prf} i) In view of Lemma \ref{2.5}, it would be enough to show that there is a balanced big Cohen-Macaulay $R$-module
$B_M$ such that $\Tor^R_1(B_M,\text{Syz}^F_i(M))=0$ for all $1\leq i\leq \dim R$. Since $\CIdim_R\ M=0$, there exists a
quasi-deformation $R\rightarrow R'\leftarrow A$ such that
$$\pd_A(R'\otimes_RM)=\pd_A(R')<\infty.$$ As $M\in \G(R)$, it possesses a complete free resolution $$F=\cdots \rightarrow F_2
\overset{d_2}{\rightarrow} F_1\overset{d_1}{\rightarrow} F_0\overset{d_0}{\rightarrow}F_{-1}\rightarrow \cdots.$$ By
\cite[Lemma 1.1]{CFH}, we can easily see that the $R'$-complex $F':=R'\otimes_RF$ is a complete free resolution of the
$R'$-module $M':=R'\otimes_RM$. In particular, $M'$ is a totally reflexive $R'$-module.
      	
Let $B_A$ be a balanced big Cohen-Macaulay $A$-module. There exists an $A$-regular sequence $\mathbf{x}:=
x_1,\ldots, x_{\ell}$, such that $R'=A/(\mathbf{x})A$. Then, it is routine to check that $B_M:=B_A/(\mathbf{x})B_A$
is a balanced
big Cohen-Macaulay $R'$-module. Namely, let $(T,\fn)$ be a local ring and $\mathbf{z}:=z_1,z_2,\ldots, z_r\in \fn$. By
\cite[Proposition A.4]{BH}, one knows that $\mathbf{z}$ is part of a system of parameters of $T$ if and only if
$\dim T/(\mathbf{z})=\dim T-r$.
Also for a finitely generated $R$-module $N$, \cite[Theorem A.11]{BH} implies that $$\dim_{R'}(R'\otimes_RN)=
\dim_RN+\dim_{R'}R'/\fm R'.$$
Putting these two facts together, we can see that a system of parameters of $R$ maps to a part of a system of parameters
for $R'$, and so
$B_M$ is also a balanced big Cohen-Macaulay $R$-module.

Let $i\in \mathbb{N}_0$ and set $N:=\Syz^{F'}_{-i-d}(M')$, where  $d:=\dim A$. Then $\Syz^{F'}_{-i}(M')=\Syz^{F'}_d(N)$.
Since $\pd_A\ M'<\infty$ and $\pd_A\ R'<\infty$, from the exact sequence $$ 0\rightarrow M' \rightarrow
R'\otimes_RF_{-1}\rightarrow
\cdots\rightarrow R'\otimes_RF_{-i-d}\rightarrow N\rightarrow 0,$$ we deduce that $\pd_AN<\infty$. By the Auslander-Buchsbaum formula,
it turns out that $\pd_AN\le d$. Let $G$ be a free resolution of the $A$-module $B_A$. By \cite[Proposition 1.1.5]{BH}, we conclude
that $G\otimes_AR'$ is a free resolution of the $R'$-module $B_M$. Let $L$ be a free resolution of the
$R$-module $\Syz^F_{-i}(M)$. As $R'$ is a faithfully flat $R$-algebra, it turns out that $R'\otimes_RL$
is a free resolution of the $R'$-module $R'\otimes_R\Syz^F_{-i}(M)$. Now, one has
\begin{align*}
\Tor^R_1(B_M,\Syz^F_{-i}(M))&\cong \Ho_{1}(B_M\otimes_RL)\\
&\cong \Ho_{1}(B_M\otimes_{R'}(R'\otimes_RL))\\
&\cong \Tor^{R'}_1\big(B_M,R'\otimes_R\Syz^F_{-i}(M)\big)\\
&\cong \Tor^{R'}_1\big(B_M,\Syz^{F'}_{-i}(M')\big)\\
&\cong \Tor^{R'}_1\big(B_M,\Syz^{F'}_{d}(N)\big)\\
&\cong \Tor^{R'}_{d+1}(B_M,N)\\
&\cong \Ho_{d+1}\big((G\otimes_A R')\otimes_{R'}N\big)\\
&\cong \Ho_{d+1}(G\otimes_AN)\\
&\cong \Tor^A_{d+1}(B_A,N)\\
&=0.      	
\end{align*}

ii) There exits an $A$-regular sequence $\mathbf{x}:=x_1,\ldots,x_{\ell}$ such that $R'=A/(\mathbf{x})A$. Let $B_A$ be a
balanced big Cohen-Macaulay $A$-module and set $B:=B_A/(\mathbf{x})B_A$. Then, as we saw above, $B$
is a balanced big Cohen-Macaulay $R$-module.

Let $M$ be a totally reflexive $R$-module and set $M':=R'\otimes_RM$.  As we saw above, by the assumption $M\in \G(R),$ it turns out
that the $R'$-module $M'$ is totally reflexive.  Then, \cite[Corollary 4.33]{AB} implies that $\Gdim_AM'<\infty.$ This implies that
$\pd_AM'<\infty$, because $A$ is assumed to be $\G$-regular. Now, by considering $B$ instead of $B_M$, we can mimic the argument given
in the part (i) to conclude that $\Tor^R_1(B,M)=0$. This completes the argument by Lemma \ref{2.5}.
\end{prf}

Now, we record the following corollary. Although it can be proved by applying Theorems \ref{2.4} and \ref{2.7}, we present a direct
argument for its first two parts. We thank an anonymous reader who suggested us this proof.

\begin{corollary}\label{2.8} Let $(R,\fm,k)$ be a local ring and let $$G=0\lo G_s\lo G_{s-1}\lo \cdots\lo G_{0}\lo 0$$
be a non-exact complex of totally reflexive $R$-modules with finite length homologies. Assume that either of the following
conditions is satisfied:
\begin{enumerate}
\item[i)]  $\CIdim_RG<\infty$; or
\item[ii)] $R$ is a quasi-specialization of a $\G$-regular local ring; or
\item[iii)] $R$ is Cohen-Macaulay.
\end{enumerate}
Then $\dim R\leq s$.
\end{corollary}

\begin{prf}  i) As $\CIdim_RG<\infty$, there is a quasi-deformation $R\rightarrow R'\leftarrow A$ such that $\CIdim_RG=\pd_A(R'\otimes_RG)-\pd_AR'$.
By localizing at a minimal prime ideal of $\fm R'$, we may and do assume that $R'/\fm R'$ is Artinian. In particular, $\dim R'=\dim R$. Since
$G$ has homology of finite length over $R$ and $R'/\fm R'$ is Artinian, \cite[Theorem A.11]{BH} yields that $R'\otimes_RG$ has homology of finite
length over $R'$ and over $A$. In particular, $\dim_A(R'\otimes_RG)\leq 0$. So, \cite[Lemma 4.2]{F}, whose proof is an application of
the existence of Big Cohen-Macaulay modules,  implies that $\dim A\leq \pd_A(R'\otimes_RG)$. Now, we have:
\begin{align*}
\dim R&=\dim R'\\
&=\dim A-\pd_AR'\\
&\leq \pd_A(R'\otimes_RG)-\pd_AR'\\
&=\CIdim_RG\\
&=\Gdim_RG\\
& \leq s.
\end{align*}

ii)  Assume that there is a quasi-deformation $R\rightarrow R'\leftarrow A$ in which $A$ is
$\G$-regular. As $A$ is $\G$-regular, one has $\Gdim_A(R'\otimes_RG)=\pd_A(R'\otimes_RG)$. Thus applying
\cite[Theorem 2.3.12]{C1} and \cite[Corollary 5.11]{C2} yield that
\begin{align*}
\pd_A(R'\otimes_RG)-\pd_AR'&=\Gdim_A(R'\otimes_RG)-\pd_AR'\\
&=\Gdim_{R'}(R'\otimes_RG)\\
&=\Gdim_RG.
\end{align*}
Therefore $\CIdim_RG<\infty$, and so the assertion follows by i).

iii) It is obvious by Theorem 2.4. Note that if $R$ has positive dimension, then by the Auslander-Bridger formula all 
totally reflexive $R$-modules have positive depth. 
\end{prf}

The following result indicates that the class of quasi-specializations of $\G$-regular local rings is quite big.

\begin{remark}\label{2.10} Let $(S,\fm_S,k)$ and  $(T,\fm_T,k)$ be two commutative Noetherian local rings with a common residue
field $k$.  Let $\pi_S:S\lo k$ and $\pi_T:T\lo k$ denote the natural epimorphisms and $$U:=\{(s,t)\in S\times T \ | \ \pi_S(s)=
\pi_T(t) \}.$$ Then $U$, by the natural pointwise multiplication, is a commutative Noetherian local ring  which is called the
fiber product
ring of $S$ and $T$.
\begin{enumerate}
\item[i)] If $U$ is not Gorenstein, then it is $\G$-regular; see \cite[Corollary 4.7]{NS}.
\item[ii)] By \cite[Fact 2.2]{NTSV},  the ring $U$ is Cohen-Macaulay if and only if $S$ and $T$ are Cohen-Macaulay and
$\dim S=\dim T\leq 1$. So, there are plenty of examples of non-Cohen-Macaulay $\G$-regular rings.
\item[iii)] Let $(R,\fm,k)$ be a $\G$-regular local ring and $x\in \fm$ an $R$-regular element. By \cite[Proposition 4.6]{T1},
the ring $R/xR$ is $\G$-regular if and only if $x\notin \fm^2$. Thus, there are plenty of examples of quasi-specializations of
$\G$-regular rings which are neither $\G$-regular nor Cohen-Macaulay.
\item[iv)] By \cite[Example 3.5 (2)]{AM}, every Golod local ring which is not a hypersurface is $\G$-regular.
\item[v)] Let $(R,\fm,k)$ be a non-Cohen-Macaulay local ring. Assume that $x_1, x_2, \ldots, x_n\in \fm$
form an $R$-regular sequence such that $\fm/(x_1, x_2, \ldots, x_n)$ is decomposable. Then $R/(x_1, x_2, \ldots, x_n)$ is a
non-Cohen-Macaulay $\G$-regular local ring; see \cite[Lemma 3.1]{Og} or \cite[Fact 3.1]{NT}.
\end{enumerate}
\end{remark}

Next, we present an example of a big Cohen-Macaulay $R$-module $B$ and a totally reflexive $R$-module $G$ such that
$\text{Soc}(B\otimes_RG)\neq 0$.

\begin{example}\label{2.12} Let $R:=k[[X,Y]]/(XY)$, where $k$ is a field. Then $R$ is a Gorenstein complete local ring of dimension
one. Denote by $x$ and $y$ the residue classes of $X$ and $Y$ in $R$ and set $B:=R/(x)$ and $G:=R/(y)$. We can easily check that $B$
and $G$ are maximal Cohen-Macaulay $R$-modules, and so both $B$ and $G$ are totally reflexive. Now, as $B\otimes _RG\cong k$, we have
$\text{Soc}(B\otimes_RG)\neq 0.$
\end{example}

\section{Question 1.2}

This section deals with Question \ref{1.2}. As, we have already mentioned Question 1.1 easily implies Question 1.2. So, we begin
with recording the following consequence of Corollary \ref{2.8} ii).

\begin{corollary}\label{3.1} Let $R$ be a quasi-specialization of a $\G$-regular local ring.  If $R$ admits  a non-zero
Cohen-Macaulay module $M$ with finite Gorenstein dimension, then $R$ is Cohen-Macaulay.
\end{corollary}

\begin{prf}   Suppose that $\mathbf{x}:=x_1,\ldots, x_t\in \fm$ is a maximal $M$-regular sequence. Then the $R$-module
$M/{(\mathbf{x})}M$ has finite length and \cite[Theorem 8.7.7]{Av} implies that $$\Gdim_R{M/{(\mathbf{x})}M}=\Gdim_RM+t.$$
Hence, we may and do assume that $M$ has finite length. Set $s:=\depth R$. By the Auslander-Bridger formula, one has
$$\Gdim_RM=\depth R-\depth_RM=s,$$ and so there is a non-exact complex $$G=0\lo G_s\lo \cdots \lo G_0\lo 0$$ of totally
reflexive $R$-modules such that $G\simeq M$. Now by Corollary \ref{2.8} ii), it turns out that $\dim R\leq s$, and so $R$ is
Cohen-Macaulay.
\end{prf}

Next, we provide an affirmative answer to Question \ref{1.2} in the case $R$ is quasi-Buchsbaum. First, we recall definition
of quasi-Buchsbaum rings.

\begin{definition}\label{3.2} Let $(R,\fm,k)$ be a local ring.
\begin{enumerate}
\item[i)] A sequence $x_1, x_2, \ldots, x_r\in \fm$ is called a {\it weak $R$-sequence} if $$\langle x_1, x_2, \ldots, x_{i-1} 
\rangle :_Rx_i=\langle x_1, x_2, \ldots, x_{i-1} \rangle :_R\fm$$ for every $i=1, 2, \ldots, r$.
\item[ii)] The local ring $R$ is said to be {\it quasi-Buchsbaum} if every system of parameters of $R$ in
$\fm^2$ is a weak $R$-sequence.
\end{enumerate}
\end{definition}

To prove the main result of this section, we need the following three lemmas.

\begin{lemma}\label{3.3} Let $(R,\fm,k)$ be a local ring. Then the following are equivalent:
\begin{enumerate}
\item[i)] There is a system of parameters of $R$ contained in $\fm^2$ which is a weak $R$-sequence.
\item[ii)] $R$ is quasi-Buchsbaum.
\item[iii)]  $\fm \text{H}_{\fm}^i(R)=0$ for all $i=0,\dots , \dim R-1$.
\end{enumerate}
\end{lemma}

\begin{prf} See \cite[Proposition 2.1]{SV}.
\end{prf}

\begin{lemma}\label{3.4} Let $(R,\fm,k)$ be a quasi-Buchsbaum local ring. For every $R$-regular sequence $x_1,x_2,\dots,
x_t\in \fm^2$, the local ring $R/\langle x_1, x_2, \ldots, x_{t} \rangle $ is quasi-Buchsbaum.
\end{lemma}

\begin{prf} Set $\overline{R}:=R/\langle x_1, x_2,\dots, x_t \rangle$ and let $-:R\to \overline{R}$ denote the natural ring
epimorphism. Set $d:=\dim R$ and let $\overline{x}_{t+1}, \overline{x}_{t+2}, \ldots, \overline{x}_d\in (\fm\overline{R})^2$
form a system of parameters of $\overline{R}$. Thus, $x_1, x_2, \dots, x_t, x_{t+1}, \dots, x_d\in \fm^2$ form a system of
parameters of $R$. As $R$ is quasi-Buchsbaum, $x_1, x_2,\dots, x_d$ is a weak $R$-sequence. In particular, $$\langle x_1, x_2,
\ldots, x_{i-1}\rangle :_Rx_i=\langle x_1, x_2, \ldots, x_{i-1} \rangle :_R\fm$$ for every $i=1, 2, \ldots, d$. Therefore,
$$\langle \overline{x}_{t+1}, \overline{x}_{t+2}, \ldots, \overline{x}_{t+i-1}\rangle :_{\overline{R}}\overline{x}_{t+i}=\langle
\overline{x}_{t+1}, \overline{x}_{t+2}, \ldots, \overline{x}_{t+i-1}\rangle:_{\overline{R}}\fm \overline{R}$$ for every $i=1, 2,
\ldots, d-t$, and so $\overline{R}$ is quasi-Buchsbaum.
\end{prf}

The next result is interesting by its own right. It improves \cite[Corollary 4.5]{As} and \cite[Lemma 5.2]{DHN} by relaxing locally free
assumption in \cite[Corollary 4.5]{As} and finite length assumption in \cite[Lemma 5.2]{DHN}.

\begin{lemma}\label{3.5} Let $(R,\fm,k)$ be a local ring for which $\Gamma_{\fm}(R)$ is a non-zero $k$-vectore space. Let $M$ be a finitely
generated $R$-module. Assume that, $\Omega_{t+1}(M)$, the $t+1$-th syzygy of $M$ has finite length for some integer $t\geq 2$. Then $M$ is free.
\end{lemma}

\begin{prf} Assume that $R$ is Artinian. Since $\Gamma_{\fm}(R)$ is a $k$-vector space, it follows that $\fm=\fm R=\fm \Gamma_{\fm}(R)=0$.
Thus $R$ is a field, and so every $R$-module is free. Hence, we may and do assume that $\dim R>0$.

Set $H:=\Gamma_{\fm}(R)$. Choose an augmented minimal free resolution $$F=\cdots\rightarrow F_2\overset{d_2}{\rightarrow}F_1\overset{d_1}{\rightarrow}
F_0\overset{d_0}{\rightarrow} M\rightarrow 0$$ of $M$. From the short exact sequence $$0\lo H\lo R\lo R/H\lo 0,  \  \   (*)$$  we obtain the following
commutative diagram with exact rows:
\vspace{.4cm}
$$\begin{CD}
0 @>>> F_t\otimes_R H @>>> F_t\otimes_R R @>>>  F_t\otimes_R R/H @>>> 0\\
 @. @Vd_t\otimes_R HVV @Vd_t\otimes_R R VV @Vd_t\otimes_R R/H VV @.\\
 0@>>> F_{t-1}\otimes_R H @>>> F_{t-1}\otimes_R R @>>> F_{t-1}\otimes_R R/H @>>>0,
\end{CD}$$\\
which, in view of the Snake Lemma, yields the exact sequence $$\Omega_{t+1}(M)\lo \ker(d_t\otimes_R R/H)\lo \coker(d_t\otimes_RH).$$ By the assumption,
$\ell_R(\Omega_{t+1}(M))<\infty$ and as $\ell_R(H)<\infty$, we get $\ell_R(\coker(d_t\otimes_RH))<\infty$. Thus $\ell_R(\ker(d_t\otimes_R R/H))<\infty$.
But $F_t\otimes_RR/H$ has positive depth, and so it can not possess any non-zero Artinian submodule. Hence $\ker(d_t\otimes_R R/H)=0$, and so $\Tor_t^R(M,R/H)=0$.

Next, the long exact Tor-sequence induced by the short exact sequence $(*)$ implies that $$\Tor_{t-1}^R(M,H)\cong \Tor_t^R(M,R/H)=0.$$ By the assumption,
$H$ is a direct sum of finitely many copies of $k$. Thus $\Tor_{t-1}^R(M,k)=0$, and so $\pd_RM\leq t-2$. As $H\neq 0$, we get $\depth R=0$. Therefore, the
Auslander-Buchsbaum formula concludes that $\pd_RM=0$.
\end{prf}

Next, we present the main result of this section.

\begin{theorem}\label{3.6} Let $(R,\fm,k)$ be a quasi-Buchsbaum local ring. Assume that $R$ admits a non-zero Cohen-Macaulay module
of finite Gorenstein dimension. Then $R$ is Cohen-Macaulay.
\end{theorem}

\begin{prf} As we saw in the proof of Corollary \ref{3.1}, we may and do assume that $R$ possesses a non-zero finite length module
$M$ of finite Gorenstein dimension. Since $M$ has finite length, it is annihilated by some power $s\geq 2$ of $\fm$. Set $t:=\depth R$.
Let $x_1,x_2,\dots, x_t\in \fm^s$ be an $R$-regular sequence and set $S:=R/\langle x_1, x_2, \ldots, x_{t} \rangle $ and $\fn:=\fm/
\langle x_1, x_2, \ldots, x_{t} \rangle $. Then $\depth S=0$ and, by Lemma \ref{3.4}, the local ring $S$ is quasi-Buchsbaum.

Once we show $\dim S=0$, the proof will be complete. On the contrary, suppose that $\dim S>0$. Then Lemma \ref{3.3} yields that
$\Gamma_{\fn}(S)$ is a non-zero $k$-vector space. Next, \cite[Corollary 4.33]{AB} implies that $\Gdim_RM=\Gdim_SM+t$. As $M$ has finite
length, by the Auslander-Bridger formula, we conclude that $\Gdim_RM=t$, and so $M$ is totally reflexive as an $S$-module.

Let $$F=\cdots\rightarrow F_2\overset{d_2}{\rightarrow}F_1\overset{d_1}{\rightarrow} F_0\overset{d_0}{\rightarrow}F_{-1}\rightarrow
\cdots,$$ be a complete free resolution of the $S$-modules $M$ and set $N:=\Syz_{-3}^F(M)$. Hence, $\Syz_{3}^F(N)=M$ and there
is an exact sequence $$0\rightarrow  M\rightarrow F_{-1}\overset{d_{-1}}{\longrightarrow} F_{-2}\overset{d_{-2}}{\longrightarrow}
F_{-3}\rightarrow N\rightarrow 0. \  \   \  \  \ \  \  (\dag)$$
Now, Lemma \ref{3.5} implies that the $S$-module $N$ is free. Hence, from $(\dag)$, one deduces that the $S$-module $M$ is also free.
This implies that the ring $S$ is Artinian, which is a contradiction.
\end{prf}

Next, we propose the following natural three questions.

\begin{question}\label{2.120} Let $(R,\fm,k)$ be a local ring. Assume that $R$ admits a non-zero totally reflexive module
of finite length. Is $R$ Artinian?
\end{question}

\begin{question}\label{2.121} Let $(R,\fm,k)$ be a local ring and $M$ a non-zero totally reflexive $R$-module.
Does $\dim_RM=\dim R$?
\end{question}

As far as we know,  the above question was first asked by Pham Hung Quy in MathOverflow.

\begin{question}\label{2.13} Let $(R,\fm,k)$ be a local ring and $M$ a totally reflexive $R$-module. Is there a big Cohen-Macaulay
$R$-module $B_M$ such that $\text{Soc}(B_M\otimes_RM)=0$?
\end{question}

\begin{remark}\label{2.141} Although, we do not know the answers to these three questions and also to the two questions in
the introduction, below we indicate the relationship between them.
\vspace{.5cm}
\begin{center}

$\xymatrix{\text{Question 3.9} \ar@{=>}[dr]^c & \text{Question 3.8} \ar@/_/@{=>}[rrr]_b \ar@/^/@{<=}[rrr]^{\text{+gCM}} &&& \text{Question 3.7}  &\\
&\text{Question 1.1}  \ar@{=>}[u]^a  \ar@{=>}[rrr]^e  &&& \text{Question 1.2}  \ar@{<=>}[u]^d}$
\end{center}
\vspace{.5cm}

\begin{enumerate}
\item[i)] Suppose that $T$ is a non-zero totally reflexive $R$-module and $\mathbf{x}$ is a system of parameters for $T$,
i.e. the sequence $\mathbf{x}$ has length $\text{dim}\ T$, and $T/\mathbf{x}T$ is of finite length. Consider the the complex,
$K_{\bullet}(\mathbf{x};T)$. At each $\mathfrak{p}\neq\mathfrak{m}$ we have $T_{\mathfrak{p}}=(\mathbf{x})T_{\mathfrak{p}}$,
thus either $(\mathbf{x})R_{\mathfrak{p}}=R_{\mathfrak{p}}$ or $T_{\mathfrak{p}}=0$ by Nakayama's Lemma. Hence, since
$\mathbf{x}$ kills the homologies $\Ho_{i}(\mathbf{x};T_{\mathfrak{p}})$, we get  $\Ho_{i}(\mathbf{x};T_{\mathfrak{p}})=0$ for
all $i$. Thus, $K_{\bullet}(\mathbf{x};T)$ has finite length homologies. If the Gorenstein version of the New Intersection
Theorem holds, then applying it to the complex $K_{\bullet}(\mathbf{x};T)$, we would have $\text{dim}\ R\le \text{dim}_R\ T$.
So, the implication (a) holds.
\item[ii)] The implication (b) is obvious.
\item[iii)] Assume that $(R,\fm,k)$ is a local ring such that $\Ass_RR\setminus \{\fm\}\subseteq \Assh_RR$
(e.g. $R$ is generalized Cohen-Macaulay; see \cite[Chapter 1, Lemma 2.2]{SV}). Then Question \ref{2.120} implies Question
\ref{2.121}. To this end, assume that the answer to
Question \ref{2.120} is affirmative and let $T$ be a non-zero totally reflexive $R$-module. Let $\fp\in \text{Assh}_RT$. Then,
we already have $\mathfrak{p}\in\text{Ass}\ R,$ because $T\cong \text{Hom}_{R}\big(\text{Hom}_{R}(T,R),R\big)$. If it is the
case that $\mathfrak{p}=\mathfrak{m}$, then $T$ would be a totally reflexive module of finite length and thus $R$ would be an
Artinian ring. On the other hand, if $\mathfrak{p}\neq\mathfrak{m}$, then $\mathfrak{p}\in\text{Assh}\ R$, and so we have
$\text{dim}\ R/\mathfrak{p}=\text{dim}\ R$.
\item[iv)] Theorem \ref{2.4} implies the implication (c).
\item[v)] The equivalence (d) follows by the argument given in the first two paragraphs of the proof of Theorem \ref{3.6}.
\item[vi)] The implication (e) is immediate by the argument which is given in the proof of Corollary \ref{3.1}.
\end{enumerate}
\end{remark}

\begin{corollary}\label{3.11} Let $(R,\fm,k)$ be a quasi-Buchsbaum ring.  Then $\dim_RM=\dim R$ for every non-zero totally reflexive
$R$-module $M$.
\end{corollary}

\begin{prf}   Follows by Theorem \ref{3.6} and Remark \ref{2.141} v) and iii). Note that any quasi-Buchsbaum ring is generalized Cohen-Macaulay.
\end{prf}


\end{document}